\documentclass[]{elsarticle}

\usepackage{amsmath, amssymb, latexsym,enumerate, mathrsfs}

\usepackage{lineno,hyperref}

\newtheorem{assumption}{\bf Assumption}
\newtheorem{remark}{\bf Remark}
\newtheorem{proposition}{\bf Proposition}
\newtheorem{corollary}{\bf Corollary}
\newtheorem{lemma}{\bf Lemma}

\def\qed{\hfill $\Box$}

\modulolinenumbers[5]

\begin{document}

\begin{frontmatter}

\title{Large deviation principle for dynamical systems coupled with diffusion-transmutation processes}

\author{Getachew K. Befekadu\corref{mycorrespondingauthor}}
\address{NRC/AFRL \& Department of Industrial System Engineering, University of Florida - REEF, 1350 N. Poquito Rd, Shalimar, FL 32579, USA}
\cortext[mycorrespondingauthor]{Corresponding author}
\ead{gbefekadu@ufl.edu}



\begin{abstract}
In this paper, we introduce a mathematical apparatus that is relevant for understanding a dynamical system with small random perturbations and coupled with the so-called transmutation process -- where the latter jumps from one mode to another, and thus modifying the dynamics of the system. In particular, we study the exit problem, i.e., an asymptotic estimate for the exit probabilities with which the corresponding processes exit from a given bounded open domain, and then formally prove a large deviation principle for the exit position joint with the type occupation times as the random perturbation vanishes. Moreover, under certain conditions, the exit place and the type of distribution at the exit time are determined and, as a consequence of this, such information also give the limit of the Dirichlet problems for the associated partial differential equation systems with a vanishing small parameter. 
\end{abstract}
\begin{keyword}
Boundary exit problem\sep diffusion process \sep large deviations\sep small random perturbations\sep transmutation process
\MSC[2010] 58J65\sep 35B20\sep 49J20\sep 60F10
\end{keyword}

\end{frontmatter}


\section{Introduction} \label{S1}
Consider the following $n$-dimensional process $x_t^{\epsilon}$ defined by
\begin{align} 
 dx_t^{\epsilon} = F(x_t^{\epsilon}, y_t^{\epsilon})dt \label{Eq1} 
\end{align}
and an $m$-dimensional diffusion process $y_t^{\epsilon}$ that obeys the following stochastic differential equation
\begin{align} 
 dy_t^{\epsilon} = f(y_t^{\epsilon})dt + \sqrt{\epsilon} \sigma(y_t^{\epsilon}) d w_t, \,\, t \in [0, T], \,\, (x_0^{\epsilon}, y_0^{\epsilon})=(x_0, y_0), \label{Eq2} 
\end{align}
where
\begin{itemize}
\item $(x_t^{\epsilon}, y_t^{\epsilon})$ jointly defined an $\mathbb{R}^{(n+m)}$-valued diffusion process,
\item the functions $F$ and $f$ are uniformly Lipschitz, with bounded first derivatives,
\item $\sigma(y)$ is a Lipschitz continuous $\mathbb{R}^{m \times m}$-valued function such that $a(y) = \sigma(y)\,\sigma^{T}(y)$ is uniformly elliptic, i.e.,
\begin{align*}
 a_{min} \vert p \vert^2 < p \cdot a(y) p < a_{max} \vert p \vert^2, \quad p, y \in \mathbb{R}^{m},
\end{align*}
for some $a_{max} > a_{min} > 0$, 
\item $w_t$ is a standard Wiener process in $\mathbb{R}^{m}$, and
\item $\epsilon$ is a small positive number that represents the level of random perturbation in the system.
\end{itemize}
Let $\tau_{G_0}^{\epsilon}$ be the first exit time for the component $x_t^{\epsilon}$ from a bounded open domain $G_0 \subset \mathbb{R}^{n}$, with smooth boundary $\partial G_0$, i.e.,
\begin{align}
\tau_{G_0}^{\epsilon} = \min \bigl\{ t > 0 \, \bigl\vert \, x_t^{\epsilon} \notin G_0 \bigr\}. \label{Eq3}
\end{align}
Here, we remark that a small random perturbation enters only in \eqref{Eq2} and is then subsequently transmitted to the other dynamical system in \eqref{Eq1}. As a result of this, the $\mathbb{R}^{(n+m)}$-valued diffusion process $(x_t^{\epsilon}, y_t^{\epsilon})$ is degenerate in the sense that the operator associated with it is a degenerate elliptic equation. Recently, the authors in \cite{BefA15a}, based on the following assumptions:
\begin{enumerate} [i)]
\item the infinitesimal generator
\begin{align*} 
\mathcal{L}_0^{\epsilon}\bigl(\cdot\bigr)(x,y) =  \frac{\epsilon}{2} \operatorname{tr}\bigl \{a(y)\bigtriangledown_y^2 \bigl(\cdot\bigr) \bigr\} + \bigl \langle \bigtriangledown_x \bigl(\cdot\bigr), F(x, y) \bigr\rangle + \bigl \langle \bigtriangledown_y \bigl(\cdot\bigr), b(y) \bigr\rangle,
\end{align*}
is hypoelliptic (e.g., see \cite{Hor67} or \cite{Ell73}); see also Remark~\ref{R2} below), and
\item $\bigl \langle F(x, y), \nu(x) \bigr\rangle > 0$, where $\nu(x)$ is a unit outward normal vector to $\partial G_{0}$ at $x \in \partial D_{0}$ and for $y \in \mathbb{R}^m$,
\end{enumerate}
have provided an asymptotic estimate for the exit probability with which the diffusion process exits from a given bounded open domain during a certain time interval. Note that the approach for such an asymptotic estimate is basically relied on the interpretation of the exit probability function as a value function for a certain stochastic control problem that is associated with the underlying dynamical system with small random perturbations.
\begin{remark}  \label{R2}
Note that, in general, the hypoellipticity assumption is related to a strong accessibility property of controllable nonlinear systems that are driven by white noise (e.g., see \cite{SusJu72} concerning the controllability of nonlinear systems, which is closely related to \cite{StrVa72} and \cite{IchKu74}; see also \cite[Section~3]{Ell73}). That is, the hypoellipticity assumption implies that the diffusion process $(x_t^{\epsilon}, y_t^{\epsilon})$ has a transition probability density with a strong Feller property.
\end{remark}

\begin{figure}[bht]
\vspace{-2mm}
\begin{center}
\hspace{20 mm}\includegraphics[width=70mm]{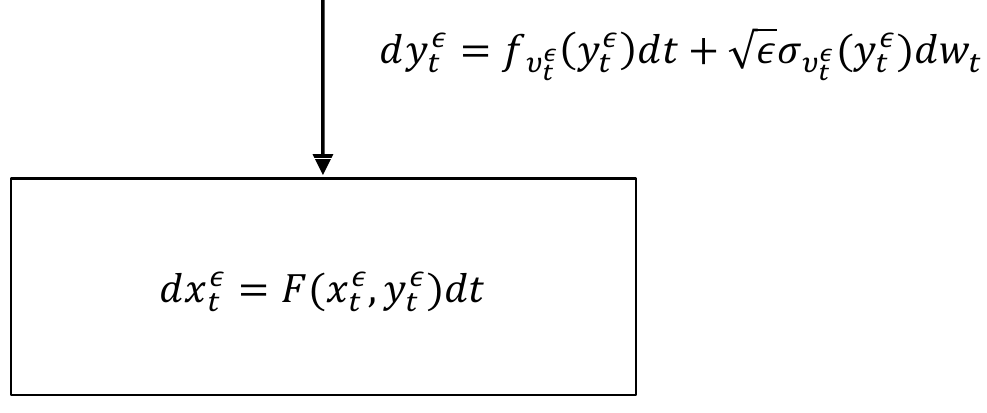}
{\scriptsize $ \begin{array}{c}
 \nu_t^{\epsilon} \,\colon\,\mathbb{P} \bigl\{ \nu_{t+\triangle}^{\epsilon} = m \, \bigl\vert \,(x_t^{\epsilon}, y_t^{\epsilon}) = (x, y), \nu_t^{\epsilon} = k \bigr\} = \dfrac{c_{km}(x,y)}{\epsilon} \triangle + o(\triangle)  \\ 
  \text{as} \,\,\, \triangle \downarrow 0  \,\,\, \text{for} \,\,\, k, m \in \{1,2, \ldots, K\} \,\,\, \text{and} \,\,\,  k \neq m
\end{array}$}
\vspace{-2 mm}
\caption{A dynamical system coupled with diffusion-transmutation processes} \label{Fig-DSC_With_DTP}
\vspace{-5 mm}
\end{center}
\end{figure}

Here, it is worth mentioning that some interesting studies on the exit probabilities for the dynamical systems with small random perturbations have been reported in literature (see, e.g., \cite{VenFre70}, \cite{FreWe84}, \cite{Kif90} and \cite{DupKu86} in the context of large deviations; see \cite{DupKu89}, \cite{Day86}, \cite{EvaIsh85}, \cite{Zab85}, \cite{Fle78}, \cite{FleTs81} and \cite{BefA15a} in connection with stochastic optimal control problems; and see \cite{Day86} or \cite{MatSc77} via asymptotic expansions approach). Note that the rationale behind our framework follows, in some sense, the settings of these papers -- where we establish a connection between the asymptotic estimates for the joint type occupation times and exit probability problems for a family of diffusion-transmutation processes and the corresponding solutions for the Dirichlet problem in a given bounded open domain with smooth boundary. Specifically, we consider the following Markov process $\bigl((x_t^{\epsilon},y_t^{\epsilon}), \nu_t^{\epsilon}\bigr)$ in the phase space $\mathbb{R}^{(n+m)} \times \{1,2, \dots, K\}$ (see Fig.~\ref{Fig-DSC_With_DTP})
\begin{align} 
\left.\begin{array}{l}
dx_t^{\epsilon} = F(x_t^{\epsilon}, y_t^{\epsilon})dt \\
dy_t^{\epsilon} = f_{\nu_t^{\epsilon}}(y_t^{\epsilon})dt + \sqrt{\epsilon} \sigma_{\nu_t^{\epsilon}}(y_t^{\epsilon}) d w_t, \quad (x_0^{\epsilon}, y_0^{\epsilon})= (x_0, y_0)
\end{array} \right \}  \label{Eq4}
\end{align}
where
\begin{itemize}
\item $(x_t^{\epsilon}, y_t^{\epsilon})$ is an $\mathbb{R}^{(n+m)}$-valued diffusion process,
\item the functions $F$ and $f_k$, $k=1,2, \ldots, K$, are uniformly Lipschitz, with bounded first derivatives,
\item $\nu_t^{\epsilon}$ is a $\{1,2, \ldots, K\}$-valued process such that 
\begin{align*}
 \mathbb{P} \Bigl\{ \nu_{t+\triangle}^{\epsilon} = m \, \bigl\vert \,(x_t^{\epsilon}, y_t^{\epsilon}) = (x, y), \nu_t^{\epsilon} = k \Bigr\} = \frac{c_{km}(x,y)}{\epsilon} \triangle + o(\triangle)  \,\,\, \text{as} \,\,\, \triangle \downarrow 0, 
\end{align*}
for $k, m \in \{1,2, \ldots, K\}$ and $k \neq m$,
\item $\sigma_k(y)$ are Lipschitz continuous $\mathbb{R}^{m \times m}$-valued functions such that  $a_k(y) = \sigma_k(y)\,\sigma_k^{T}(y)$ are uniformly elliptic, that is,
\begin{align*}
 a_{min} \vert p \vert^2 < p \cdot a_k(y) p < a_{max} \vert p \vert^2, \quad p, y \in \mathbb{R}^{m}, k \in \{1,2, \ldots, K\},
\end{align*}
for some $a_{max} > a_{min} > 0$, 
\item $w_t$ is a standard Wiener process in $\mathbb{R}^{m}$, and
\item $\epsilon$ is a small positive number that represents the level of random perturbation in the system.
\end{itemize}

\begin{remark}
Note that the coupling of such a diffusion process with that of the so-called transmutation process in \eqref{Eq4} allows us to model random jumps or switchings from one state or mode to another, and thus modifying the dynamics of the systems.
\end{remark}

Define $z_t^{\epsilon}=(x_t^{\epsilon}, y_t^{\epsilon})$, then, with minor abuse of notation, we can rewrite equation~\eqref{Eq4} as follows
\begin{align} 
dz_t^{\epsilon} = F_{\nu_t^{\epsilon}}(z_t^{\epsilon})dt+ \sqrt{\epsilon} \hat{\sigma}_{\nu_t^{\epsilon}}(z_t^{\epsilon}) d w_t, \quad z_0^{\epsilon}=(x_0, y_0), \label{Eq51}
\end{align}
where
\begin{align*} 
\left.\begin{array}{l}
F_k = \operatorname{blockdiag}\bigl\{F,\, f_k\bigr\} \\
\hat{\sigma}_k = \bigl[0_{n \times n}, \sigma_k^T\bigr]^T
\end{array} \right \}
\end{align*}
for $k = 1, 2, \ldots, K$.

Here, we also assume that the transmutation coefficients $c_{km}(z)$, for $(z) \in \mathbb{R}^{(n+m)}$, are positive and Lipschitz continuous. Moreover, under these conditions (e.g., see \cite{EizF90}), there exists a unique vector $\bar{\omega}(z) = \bigl(\omega_1(z), \omega_2(z), \ldots, \omega_K(z)\bigr)$ such that
\begin{align*}
 \omega_k(z) > 0, \quad \sum\nolimits_{k = 1}^{K} \omega_k(z) = 1 \quad \text{and} \quad \bar{\omega}(z) \hat{C}(z) = 0,
\end{align*}
where $\hat{C}(z) = \bigl(\hat{C}_{km}(z)\bigr)$ is an $K \times K$ matrix and
\begin{align*}
\left \{ \begin{array}{l}
 \hat{C}_{km}(z) = c_{km}(z) \quad\qquad \quad ~\text{for} \quad k \neq m,\\
 \hat{C}_{kk}(z) = - \sum\nolimits_{j: j \neq k} c_{kj}(z) \quad \text{for}  \quad k = m.
\end{array} \right.
\end{align*}

Denote by $\mathbb{P}_{z_0, k}^{\epsilon}$ the probability measures in the space of trajectories of the process $(z_t^{\epsilon}, \nu_t^{\epsilon})$ and by $\mathbb{E}_{z_0, k}^{\epsilon}$  the associated expectation. Define the occupation time $r_t^{\epsilon}$ for the component $\nu_t^{\epsilon}$ as 
\begin{align*}
 r_t^{\epsilon} = \biggl(\int_0^t \chi_1\bigl(\nu_s^{\epsilon}\bigr) ds, \int_0^t \chi_2\bigl(\nu_s^{\epsilon}\bigr) ds, \ldots, \int_0^t \chi_K\bigl(\nu_s^{\epsilon}\bigr) ds \biggr),
\end{align*}
where $\chi_k$ is the indicator function of the singleton set $\{k\}$. Then, we specifically study the process $(z_t^{\epsilon}, \nu_t^{\epsilon})$ and the occupation time $r_t^{\epsilon}$; and we further look on $z_t^{\epsilon}$ as a result of small random perturbations of the average system
\begin{align*}
 \dot{z}(t) &= \sum\nolimits_{k=1}^K \omega_k(z(t)) F_k \bigl(z(t)\bigr)  \notag \\
                        &\triangleq \bar{F} \bigl(z(t)\bigr), \quad \quad z(0) = z_0 \in \mathbb{R}^{(n+m)}
\end{align*}
that allows us to prove large deviation results for the joint type occupation times and positions as $\epsilon \rightarrow 0$ and study the exit probabilities for such a family of processes.

On the other hand, let $G \subset \mathbb{R}^{(n+m)}$ be a bounded open domain with smooth boundary $\partial G$. The infinitesimal generator $\mathcal{L}^{\epsilon}$ of the process $(z_s^{\epsilon}, \nu_s^{\epsilon})$ acting on smooth functions (smooth in $z \in \mathbb{R}^{(n+m)}$) is given by
\begin{align}
 \mathcal{L}^{\epsilon}\upsilon_k(z)= \mathcal{L}_k^{\epsilon}\upsilon_k(z) + \frac{1}{\epsilon} \sum\nolimits_{j=1}^K c_{kj}(z) \bigl[\upsilon_j(z) - \upsilon_k(z)\bigr], \label{Eq4b}
\end{align}
where
\begin{align}
\mathcal{L}_{k}^{\epsilon} \upsilon_k(z) = \Bigl \langle \bigtriangledown_z \upsilon_k(z), F_k(z) \Bigr\rangle + \frac{\epsilon}{2} \operatorname{tr}\bigl \{\hat{a}_k(z)\bigtriangledown_z^2 \upsilon_k(z) \bigr\}, \label{Eq4c}
\end{align}
and $\hat{a}_k(z) \triangleq \bigl(\hat{a}_k^{ij}(z)\bigr) = \hat{\sigma}_k(z)\,\hat{\sigma}_k^T(z)$ for $k =1, 2, \ldots, K$.  Moreover, throughout the paper, we assume that the infinitesimal generator in \eqref{Eq4b} is hypoelliptic for each $k =1, 2, \ldots, K$ (cf. Remark~\ref{R2}). Note that the process $(z_s^{\epsilon}, \nu_s^{\epsilon})$ is closely connected with the following Dirichlet problem for a linear reaction-diffusion system, which also satisfies the maximum principle (e.g., see \cite[Chapter~3, Section~8]{ProW84} for the application of maximum principle for classical Dirichlet problems),
\begin{align}
\left \{ \begin{array}{l}
  \mathcal{L}_k^{\epsilon} \upsilon_k^{\epsilon}(z) + \frac{1}{\epsilon} \sum\nolimits_{j=1}^K c_{kj}(z) \bigl[\upsilon_j^{\epsilon}(z) - \upsilon_k^{\epsilon}(z)\bigr] = 0, \quad z \in G,\\
 \upsilon_k^{\epsilon}(z)\vert_{\partial G} = g_k(z), \quad k=1,2, \ldots, K,
\end{array} \right. \label{Eq5}
\end{align}
where, we can study the limiting behavior for the solution of the Dirichlet problem in \eqref{Eq5} as the small random perturbation vanishes, i.e., $\epsilon \rightarrow 0$. Here, we remark that the interplay between the small diffusion and the jumps $\nu$-component leads to the situation, where $g_k(z)$, for $k=1, 2, \ldots, K$, will influence the $\lim_{\epsilon \downarrow 0} \upsilon_k^{\epsilon}(z)$. 

In this paper, we also investigate the limiting behavior for the solution of the Dirichlet problem in \eqref{Eq5} in two steps: (i) the first step is related with the exit problem for the component $z_t^{\epsilon}$ from the domain $G$, where such an exit problem can be addressed by determining the action functional for the family of processes $z_t^{\epsilon}$ as $\epsilon \rightarrow 0$, and (ii) the second step is related with determining the position of the fast component $\nu_t^{\epsilon}$ at the random time $\tau_{G}^{\epsilon} = \bigl\{t > 0 \,\bigl\vert \, z_t^{\epsilon} \notin \partial G \bigr\}$.

The rest of the paper is organized as follows. In Section~\ref{S2}, using the basic remarks made in Sections~\ref{S1}, we briefly discuss the action functional for a dynamical system with small random perturbations and coupled with the so-called transmutation process. In this section, we also discuss the Dirichlet problem corresponding to a linear reaction-diffusion system w.r.t. a given bounded open domain. In Section~\ref{S3}, we provide our results on the asymptotic estimates for the joint type occupation times and exit probabilities for a family of diffusion-transmutation processes and the corresponding solutions for the Dirichlet problem in a given bounded open domain with smooth boundary. 

\section{Action functional for the family $(z_t^{\epsilon}, r_t^{\epsilon})$} \label{S2}
In this section, we provide some preliminary results that are concerned with the action functional for the family of processes $(z_t^{\epsilon}, r_t^{\epsilon})$ as $\epsilon$ tends to zero. Before stating these results, we need some notations. Let $\lambda(z, p, \alpha)$ be the principal eigenvalue of the matrix $\bigl(H_{km}(z, p,\alpha)\bigr)$, $z, p \in \mathbb{R}^{(n+m)}$, $\alpha = (\alpha_1, \alpha_2, \ldots, \alpha_K) \in \mathbb{R}^K$
\begin{align}
H_{km}(z, p,\alpha) = \left \{ \begin{array}{l}
  \hat{C}_{km}, \hspace{2.07in} \text{if} \quad m \neq k,\\
  \bigl[p \cdot \hat{a}_k(z) p/2 + p \cdot F_k(z) + \alpha_k \bigr] + \hat{C}_{kk}, \quad \text{if} \quad m = k
\end{array} \right. \label{Eq6}
\end{align}
where $\hat{a}_k(z) = \hat{\sigma}_k(z)\,\hat{\sigma}_k^T(z)$ for $k =1, 2, \ldots, K$.

Note that $\lambda(z, p, \alpha)$ is convex in $(p, \alpha)$ and its Legendre transform  in $(p, \alpha)$ is given by
\begin{align}
\eta(z, q, \beta) =\sup_{p \in \mathbb{R}^{(n+m)},\, \alpha \in \mathbb{R}^K}  \bigl[ q \cdot p + \beta \cdot \alpha - \lambda(z, p, \alpha) \bigr], \quad z, q \in \mathbb{R}^{(n+m)}, \,\,\beta \in \mathbb{R}^K. \label{Eq7}
\end{align}

Let $C(\mathbb{R}^{(n+m)})$ be the space of continuous functions: $[0, T] \rightarrow \mathbb{R}^{(n+m)}$ and 
\begin{align}
C_{+}(\mathbb{R}^K) = \Bigl\{\mu=(\mu_1, \mu_2, \ldots, \mu_K) \,\bigl \vert \,\mu \in C(\mathbb{R}^K), \quad  \mu_k(0)=0, \,\,  \,\, 1 \le k \le K, \quad  & \notag \\
  \mu_k(t)  \,\, \text{is non-decreasing and} \,\, \sum\nolimits_{k=1}^K \mu_k(t) = t, \,\, t \in [0, T] \Bigr\}.& \label{Eq8}
\end{align}

Let $T > 0$ be fixed and define
\begin{align}
S_{0T}(\varphi, \mu) = \left \{ \begin{array}{l}
  \int_0^T \eta\bigl(\varphi(s), \dot{\varphi}(s), \dot{\mu}(s)\bigr) ds, \quad \text{if} \,\, \varphi \in C(\mathbb{R}^{(n+m)}) \,\, \text{and} \,\, \mu \in C_{+}(\mathbb{R}^K) \\
  \hspace{1.8 in} \text{are absolutely continuous (a.c.)}\\
 +\infty \hspace{1.65 in}\text{otherwise}.
\end{array} \right. \label{Eq9}
\end{align}

Suppose that the diffusion and transmutation coefficients satisfy the Lipschitz continuous and positive-Lipschitz continuous conditions, respectively. Then, we have the following result.
\begin{proposition} \label{P1}
The functional $\epsilon^{-1} S_{0T}$ is the action functional for the family of processes $(x_t^{\epsilon}, r_t^{\epsilon})$ as $\epsilon \rightarrow 0$ in the uniform topology. The action $S_{0T}$ is nonnegative and equal to zero only when $\dot{\varphi}(t) = \bar{F} \bigl(\varphi(t)\bigr) = \sum_{k+1}^K \omega_k(\varphi(t)) F_k \bigl(\varphi(t)\bigr)$ and $\dot{\mu}(t) = \bar{\omega}(\varphi(t))$, for $t \in [0, T]$.
\end{proposition}

Let us denote by $\Psi_t^{z_0}$ the integral curve of the vector field $\dot{z}(t)=\bar{F} \bigl(z(t)\bigr)$ starting from the point $z(0)=z_0$ (i.e., $\dot{\Psi}_t^{z_0} = \bar{F} \bigl( \Psi_t^{z_0}\bigr)$, with $\Psi_0^{z_0} = z_0$). Then, define
\begin{align*}
\rho(z, q) =\sup_{p \in \mathbb{R}^{(n+m)}}  \bigl[ q \cdot p - \lambda(z, p, 0) \bigr], \quad z, q \in \mathbb{R}^{(n+m)}
\end{align*}
and
\begin{align*}
I_{0T}(\varphi) = \left \{ \begin{array}{l}
  \int_0^T \rho\bigl(\varphi(s), \dot{\varphi}(s)\bigr) ds, \quad \quad \text{if} \,\, \varphi \in C(\mathbb{R}^{(n+m)}) \,\, \text{is a.c.}\\
 +\infty \hspace{1.2 in}\text{otherwise}.
\end{array} \right.
\end{align*}
Taking into account the involution property of the Legendre transform, then we have the following
\begin{align*}
\eta(z, q, \beta) &= -\sup_{\beta \in \mathbb{R}^K}  \bigl[ - 0 \cdot \beta - \eta(z, q, \beta)  \bigr] \\
                         &= \sup_{p \in \mathbb{R}^{(n+m)}}  \bigl[ q \cdot p - \lambda(z, p, 0) \bigr]\\
                         &= \rho(x, q). 
\end{align*}

Next, we have the following result which is a direct consequence of the contraction principle (see also \cite[Chapter~5, pp.~117--124]{FreWe84}).
\begin{corollary}  \label{C1}
The functional $\epsilon^{-1} I_{0T}$ is the action functional for the family $z_t^{\epsilon}$ as $\epsilon \rightarrow 0$ in the uniform topology. The action $I_{0T}$ is equal to zero only when $\varphi_t = \Psi_t^{z_0}$ and $\Psi_0^{z_0}=z_0$.
\end{corollary}
Let $\bar{n}(y)$ be a unit vector normal to $G$ at $y \in \partial G$. Furthermore, we assume that the average system $\bar{F}(z)=\sum\nolimits_{k=1}^K \omega_k(z) F_k(z)$ satisfies following large deviation condition.

\begin{assumption} [Large deviation condition] \label{AS1}
The vector field $\bar{F}(z)$ points inward from the boundary $\partial G$, i.e., $\bigl\langle \bar{F}(z), \bar{n}(z)\bigr\rangle < 0$ for any $z \in \partial G$. The vector field $\bar{F}(z)$ has a unique stationary point at $\bar{z}_0 \in G$. Moreover, the function
\begin{align*}
V(z) = \inf \Bigl\{ I_{0T}(\varphi) \, \bigl\vert \, \varphi(0) = z_0, \,\, T > 0, \,\, \varphi(T) = z \,\,\, \text{for} \,\,\, z \in \partial G \Bigr\}
\end{align*}
attains its unique minimum at $\bar{z}_0 \in \partial G$, i.e., $V(\bar{z}_0) < V(z)$ for any $z \in \partial G$.
\end{assumption}

\begin{assumption} \label{AS2}
There exists $k_0$, with $k_0 \in \{1,2, \ldots, K\}$, such that at the point $\bar{z}_0 \in \partial G$, defined above in Assumption~\ref{AS1}, then the following generic inequalities hold
\begin{align*}
 \bigl\langle F_{k_0}(\bar{z}_0), \bar{n}(\bar{z}_0)\bigr\rangle > \bigl\langle F_k (\bar{z}_0), \bar{n}(\bar{z}_0)\bigr\rangle, \,\, k_0 \in \{1,2, \ldots, K\}, \,\, k \neq k_0.
\end{align*}
Moreover, we also assume that the infinitesimal generator in \eqref{Eq4b} is hypoelliptic for each $k \in \{1, 2, \ldots, K\}$.
\end{assumption}

Let $\tau_{G}^{\epsilon}$ be the first exit time for the component $z_t^{\epsilon}$ from $G \times \{1,2, \ldots, K\}$, i.e.,
\begin{align}
\tau_{G}^{\epsilon} = \min \Bigl\{ t > 0 \, \bigl\vert \, z_t^{\epsilon} \notin G \Bigr\}. \label{Eq10}
\end{align}
Then, in the following sections, we study the limiting distribution of $(z_{\tau_{G}^{\epsilon}}^{\epsilon}, \nu_{\tau_{G}^{\epsilon}}^{\epsilon})$ as $\epsilon \rightarrow 0$. This distribution also defines the limit for the solution of the Dirichlet problem in \eqref{Eq5} as $\epsilon \rightarrow 0$.

\section{Main results} \label{S3}
In this section, we present our main results that establish connection between the exit probability problem for $(z_t^{\epsilon}, \nu_t^{\epsilon})$ from $G \times \{1, 2, \ldots, K\}$ and that of the limiting behavior for the solutions of the Dirichlet problem in \eqref{Eq5}. Note that if Assumption~\ref{AS1} holds true (i.e., the large deviation condition), then the exit problem for the component $z_t^{\epsilon}$ from the given bounded open domain $G$ is equivalent to determining the action functional for the family of processes $z_t^{\epsilon}$ as $\epsilon \rightarrow 0$ (cf. Proposition~\ref{P1}) and the exact exit position for the fast component $\nu_t^{\epsilon}$ at the random time $\tau_{G}^{\epsilon} = \bigl\{t > 0 \,\bigl\vert \, z_t^{\epsilon} \notin \partial G \bigr\}$.

Then, we have our first result concerning the asymptotic estimates for the joint type occupation times and the exit positions.
 \begin{proposition} \label{P2}
Let the diffusion matrices $a_k(x)$ and transmutation coefficients $c_{km}(z)$ be Lipschitz continuous and let $a_k(z)$ be uniformly elliptic and $c_{km}(z) > 0$ for $z \in G \cup \partial G$, $k, m \in \{1,2, \ldots, K\}$, with $k \neq m$. If Assumption~\ref{AS1} holds true, i.e., the large deviation condition. Then, we have
\begin{align}
  \lim_{\epsilon \rightarrow 0} \mathbb{P}_{x_0, k}^{\epsilon} \bigl \{ \vert z_{\tau_{G}^{\epsilon}}^{\epsilon} - \bar{z}_0 \vert > \delta \bigr\} = 0 \label{Eq11}
\end{align}
for any $\delta > 0$, $1 \le k \le K$, uniformly in $\bar{z}_0 \in \Omega$ for any compact $\Omega \subset G$. Furthermore, if Assumption~\ref{AS2} is satisfied, then we have the following
\begin{align}
  \lim_{\epsilon \rightarrow 0} \mathbb{P}_{z_0, k}^{\epsilon} \bigl \{\nu_{\tau_{G}^{\epsilon}}^{\epsilon} = k_0 \, \vert \, \tau_{G}^{\epsilon} < \infty \bigr\} = 1 \label{Eq12}
\end{align}
for  $1 \le k \le K$ and $\bar{z}_0 \in \Omega \subset G$.
\end{proposition}

Before attempting to prove Proposition~\ref{P2}, let us consider the following non-degenerate diffusion process $(x_t^{\epsilon,\hat{\delta}}, y_t^{\epsilon,\hat{\delta}})$ satisfying
\begin{align}
\left.\begin{array}{l}
d x_t^{\epsilon,\hat{\delta}} = F\bigl(x_t^{\epsilon,\hat{\delta}}, y_t^{\epsilon,\hat{\delta}}\bigr) dt + \sqrt{\hat{\delta}} dv_t \\
d y_t^{\epsilon,\hat{\delta}} = f_{\nu_t^{\epsilon,\hat{\delta}}}\bigl(y_t^{\epsilon,\hat{\delta}}\bigr) dt + \sqrt{\epsilon} \sigma_{\nu_t^{\epsilon,\hat{\delta}}} \bigl(y_t^{\epsilon,\hat{\delta}}\bigr) dw_t
\end{array} \right\} \label{Eq10x}
\end{align}
with an initial condition $\bigl(x_0^{\epsilon,\hat{\delta}}, y_0^{\epsilon,\hat{\delta}}\bigr) = \bigl(x_0, y_0\bigr)$ and $v_t$ (with $v_0=0$) is an $n$-dimensional standard Wiener process and independent to $w_t$. Moreover, for the case $\epsilon^{-1}\hat{\delta} \rightarrow 0$ as $\hat{\delta} \rightarrow 0$, then we can assume that $\nu_t^{\epsilon,\hat{\delta}}$ is a $\{1,2, \ldots, K\}$-valued process satisfies the following
\begin{align*}
 \mathbb{P} \Bigl\{ \nu_{t+\triangle}^{\epsilon,\hat{\delta}} = m \, \bigl\vert \,(x_t^{\epsilon, \hat{\delta}}, y_t^{\epsilon, \hat{\delta}}) = (x, y), \nu_t^{\epsilon,\hat{\delta}} = k \Bigr\} = \frac{c_{km}(x,y)}{\epsilon} \triangle + o(\triangle)  \,\,\ \text{as} \,\,\, \triangle \downarrow 0, 
\end{align*}
for $k, m \in \{1,2, \ldots, K\}$ and $k \neq m$.

Define $z_t^{\epsilon,\hat{\delta}}=(x_t^{\epsilon,\hat{\delta}}, y_t^{\epsilon,\hat{\delta}})$, then we can rewrite the system equations in \eqref{Eq10x} as follows
\begin{align} 
dz_t^{\epsilon,\hat{\delta}} = \tilde{F}_{\nu_t^{\epsilon,\hat{\delta}}}(z_t^{\epsilon,\hat{\delta}})dt+ \sqrt{\epsilon} \tilde{\sigma}_{\nu_t^{\epsilon,\hat{\delta}}}(z_t^{\epsilon}) d \tilde{w}_t, \quad z_0^{\epsilon}=(x_0, y_0), \label{Eq51}
\end{align}
where $\tilde{w}_t = \left[v_t^T, \, w_t^T \right]^T$ and 
\begin{align*} 
\left.\begin{array}{l}
\tilde{F}_k = \operatorname{blockdiag}\bigl\{F,\, f_k\bigr\} \\
\tilde{\sigma}_k = \left[\sqrt{(\hat{\delta}/\epsilon)}\, I_{n \times n},\, \sigma_k^T\right]^T
\end{array} \right \}
\end{align*}
for $k = 1, 2, \ldots, K$.

Let $\tau_{G}^{\epsilon,\hat{\delta}}$ be the first exit time for the diffusion process $z_t^{\epsilon,\hat{\delta}}$ from the domain $G$. Moreover, define the following
\begin{align*}
\begin{array}{c}
   \theta = \tau_{G}^{\epsilon} \wedge T, \quad\quad  \theta^{\hat{\delta}} = \tau_{G}^{\epsilon,\hat{\delta}} \wedge T, \\
   \bigl\Vert z^{\epsilon,\hat{\delta}} - z^{\epsilon} \bigr\Vert_t = \sup\limits_{s \le r \le t} \Bigl\vert z_r^{\epsilon,\hat{\delta}} - z_r^{\epsilon} \Bigr\vert \,\,\, \text{for} \,\,\, s, t \in [0, T].
   \end{array}
\end{align*}
Then, we need the following lemma, which is useful for the development of our main results.
\begin{lemma} (cf. \cite[Lemma~2.5]{BefA15a}) \label{L1}
Suppose that $\epsilon > 0$ is fixed. Then, for any initial point $z_s^{\epsilon} \in G$, with $t > s$, the following statements hold true
\begin{enumerate} [(i)]
\item $\bigl\Vert z^{\epsilon,\hat{\delta}} - z^{\epsilon} \bigr\Vert_t  \rightarrow 0$, 
\item $\theta^{\hat{\delta}} \rightarrow \theta$ \, and \, $z_{\theta^{\hat{\delta}}}^{\epsilon,\hat{\delta}} \rightarrow z_{\theta}^{\epsilon}$,
\end{enumerate}
as $\hat{\delta} \rightarrow 0$, almost surely.
\end{lemma}

Here, we remark that, for the case $\epsilon^{-1}\hat{\delta} \rightarrow 0$ as $\hat{\delta} \rightarrow 0$, the solutions for the Dirichlet problem that correspond to the non-degenerate diffusion process $(x_t^{\epsilon,\hat{\delta}}, y_t^{\epsilon,\hat{\delta}})$ continuity converge to the solutions of \eqref{Eq5} as the limiting case, when $\hat{\delta} \rightarrow 0$ (cf. Lemma~\ref{L1}, Parts~(i) and (ii)); and see also the Lebesgue's dominated convergence theorem (see \cite[Chapter~4]{Roy88}).
 
Next, let us establish the following results (i.e., Proposition~\ref{P3} and Proposition~\ref{P4}) that are useful for proving the Proposition~\ref{P2}.
\begin{proposition} \label{P3}
Suppose that the drift, diffusion and transmutation coefficients are independent of the position variable $z$ (i.e., $F_k(z)$, $\sigma_k(z)$ and $c_{km}(z)$, $k, m \in \{1, 2, \ldots, K\}$, are constants). Then, the statement in Proposition~\ref{P1} holds true. 
\end{proposition}

{\em Proof}: 
Let us first simplify some notation: $\lambda(p, \alpha)$ and $\eta(q, \beta)$ for $\lambda(z, p, \alpha)$ and $\eta(z, q, \beta)$, respectively. Then, define
\begin{align*}
L(\beta) =\sup_{\alpha}  \bigl[\beta \cdot \alpha - \lambda(0, \alpha) \bigr], \quad \beta \in \mathbb{R}^K.
\end{align*}
The action functional for the family of processes $r_t^{\epsilon}$, $\epsilon \rightarrow 0$ (assuming that $\epsilon^{-1}\hat{\delta} \rightarrow 0$ as $\hat{\delta} \rightarrow 0$), is known to be ${\epsilon}^{-1}R_{0t}\colon C_{+}(\mathbb{R}^K) \rightarrow [0, \infty]$ (cf. \cite[Chapter~7]{FreWe84}), where 
\begin{align*}
R_{0t}(\mu) = \left \{ \begin{array}{l}
  \int_0^t L\bigl(\dot{\mu}(s)\bigr) ds, \quad\quad \text{if} \,\, \mu \in C_{+}(\mathbb{R}^K)\,\, \text{is a.c.}\\
 +\infty \hspace{1.0 in}\text{otherwise}.
\end{array} \right.
\end{align*}
Let $y_t^{\epsilon,\hat{\delta}, k}$, for $k \in \{1,2, \ldots, K\}$, with $y_0^{\epsilon,\hat{\delta}, k}=0$, be independent diffusion processes corresponding to the generators
$\bigl \langle \bigtriangledown_z(\cdot), \tilde{F}_k \bigr\rangle + \frac{\epsilon}{2} \operatorname{tr}\bigl \{\tilde{a}_k\bigtriangledown_z^2(\cdot) \bigr\}$, with $\tilde{a}_k= \tilde{\sigma}_k \tilde{\sigma}_k^T$ and $k \in \{1,2, \ldots, K\}$. Then, the action functional for the family of processes $y_t^{\epsilon,\hat{\delta}, k}$ is known (cf. \cite[Chapter~3]{FreWe84}), i.e., the action functional $\epsilon^{-1}R_{0t}^k(\varphi)\colon C_{+}(\mathbb{R}^n) \rightarrow [0, \infty]$, with
\begin{align*}
R_{0t}^{k}(\varphi(t)) = \left \{ \begin{array}{l}
  \frac{1}{2} \int_{0}^{t}\bigl\Vert\dot{\varphi}(s) - \tilde{F}_k \bigr\Vert_{[\tilde{a}_k]^{-1}}^2ds \quad\quad \text{if} \,\, \varphi \,\,\text{is a.c.}\\
 +\infty \hspace{1.67 in}\text{otherwise},
\end{array} \right.
\end{align*}
where $\Vert \cdot \Vert_{[\tilde{a}_k]^{-1}}$ denotes the Riemannian norm of a tangent vector.

Note that if  $r_t^{\epsilon, \hat{\delta}} = (r_{1, t}^{\epsilon,\hat{\delta}}, r_{2, t}^{\epsilon,\hat{\delta}}, \ldots r_{K, t}^{\epsilon,\hat{\delta}})$ is independent of $y_t^{\epsilon, \hat{\delta}, k}$ for $k \in \{1,2, \ldots, K\}$, then it is clear that the process $z_t^{\epsilon, \hat{\delta}}$ can be realized as
\begin{align*}
 z_t^{\epsilon,\hat{\delta}} = \sum\nolimits_{k=1}^K  r^{\epsilon, \hat{\delta}, k} \bigl(r_{k, t}^{\epsilon,\hat{\delta}}\bigr), \quad \text{with} \quad y^{\epsilon, \hat{\delta}, k} \bigl(s\bigr) = y_s^{\epsilon,\hat{\delta}, k}.
\end{align*}
Moreover, the following transformation $S$
\begin{align*}
\bigl(y_t^{\epsilon,\hat{\delta},1}, y_t^{\epsilon,\hat{\delta}, 2}, \ldots, y_t^{\epsilon, \hat{\delta}, K}, r_{ t}^{\epsilon,\hat{\delta}}\bigr) \xrightarrow{S} \bigl(z_{ t}^{\epsilon,\hat{\delta}}, r_{ t}^{\epsilon,\hat{\delta}} \bigr)
\end{align*}
is continuous in the uniform topology on any finite interval $[0, T]$ (cf. Lemma~\ref{L1}). Then, for a fixed $T>0$, it follows from \cite[Theorem~3.3.1]{FreWe84} that 
\begin{align*}
\tilde{S}_{0T}(\varphi, \mu) = \inf\biggl\{\int_0^T & L\bigl(\dot{\mu}(s)\bigr) ds + \frac{1}{2}  \sum\nolimits_{k=1}^K \int_{0}^{\mu_k(t)}\bigl\Vert\dot{\varphi}_k(s) - \tilde{F}_k \bigr\Vert_{[\tilde{a}_k]^{-1}}^2ds, \\ 
&  \text{if} \,\, \varphi_k \,\, \text{are a.c. and} \,\,\, S\bigl(\varphi_1, \varphi_2, \ldots, \varphi_K, \mu \bigr)=\bigl(\varphi, \mu \bigr) \biggr\}
\end{align*}
is the action functional for the family of processes $\bigl(z_{ t}^{\epsilon,\hat{\delta}}, r_{ t}^{\epsilon,\hat{\delta}} \bigr)$ as $\epsilon \rightarrow 0$. Note that $\tilde{S}_{0T}(\varphi, \mu)$ is equal to zero only when $\dot{\mu}(t) = \bar{\omega}$ and $\dot{\varphi}(t) = \sum\nolimits_{k=1}^K \omega_k \tilde{F}_k$ for $t \in [0, T]$.   

It remains to show that
\begin{align*}
\tilde{S}_{0T}(\varphi, \mu) = S_{0T}(\varphi, \mu) \quad \text{for a.c.} \,\,  \varphi \,\,  \text{and} \,\, \mu. 
\end{align*}
To prove this, we will first use the Jensen's inequality to get the following
\begin{align}
&\tilde{S}_{0T}(\varphi, \mu) \ge   \inf_{\varphi_1(\mu_1(T)), \varphi_2(\mu_2(T)), \ldots, \varphi_K(\mu_K(T))} \biggl\{T L\bigl(\mu(T)/T\bigr) ds \notag \\
& \quad + \sum\nolimits_{k=1}^K \frac{1}{2\mu_k(T)}\Bigl\Vert\varphi_k(\mu_k(T)) - \mu_k(T)\tilde{F}_k \Bigr\Vert_{[\tilde{a}_k]^{-1}}^2 \, \biggl \vert \, \varphi(T) = \sum\nolimits_{k=1}^K \varphi_k(\mu_k(T)) \biggr\}. \label{Eq14}
\end{align}
We claim that the right-hand side is equal to
\begin{align*}
T\lambda(p, \alpha)=& \sup_{q \in \mathbb{R}^{(n+m)},\, \beta \in \mathbb{R}^K}  \biggl\{ p \cdot q + \alpha \cdot \beta \\
 & \quad- \inf_{q_1, q_2, \ldots, q_K\colon \sum\nolimits_{k=1}^K q_k=q} \biggl[ T L(\beta/T)  + \sum\nolimits_{k=1}^K \frac{1}{2 \beta_k}\Bigl\Vert q_k - \beta_k \tilde{F}_k \Bigr\Vert_{[a_k]^{-1}}^2 \biggr] \biggr\}.
\end{align*}
This is done as follows: the right-hand side of the last equality is equal to
\begin{align*}
 \sup_{\beta} &  \biggl\{ \sup_{q_1, q_2, \ldots, q_k} \biggl[ p (q_1 + q_2 + \ldots + q_K)  - \sum\nolimits_{k=1}^K \frac{1}{2 \beta_k}\Bigl\Vert q_k - \beta_k \tilde{F}_k \Bigr\Vert_{[a_k]^{-1}}^2 \biggr] + \alpha \cdot \beta - T L(\beta/T) \biggr\}\\
 & = \sup_{\beta}  \biggl\{ \biggl(\sum\nolimits_{k=1}^K \beta_k \bigl(p \cdot a_k p/2 + p \cdot \tilde{F}_k\bigr) \biggr) + \alpha \cdot \beta - T L(\beta/T) \biggr\}\\
 & = T\lambda(p, \alpha).
\end{align*}
So far, we have shown that
\begin{align*}
\tilde{S}_{0T}(\varphi, \mu) \ge T\eta(\varphi(T)/T, \mu(T)/T).
\end{align*}
On other hand, if we consider the following restriction $\varphi(t_j) = \sum\nolimits_{k=1}^K \varphi(\mu(t_j))$ in \eqref{Eq14} for $1 \le j \le N$, with $t_0=0 < t_1 < \cdots < t_N =T$ and repeated the same argument for each of the intervals $[t_j, t_{j+1}]$. Then, we will have
\begin{align}
\tilde{S}_{0T}(\varphi, \mu) \ge \sum\nolimits_{j=1}^N (t_j - t_{j-1}) \eta\biggl(\frac{\varphi(t_j) - \varphi(t_{j-1})}{t_j - t_{j-1}},  \frac{\mu(t_j) - \mu(t_{j-1})}{t_j - t_{j-1}} \biggr). \label{Eq15}
\end{align}
Since the partition $\{t_0, t_1, \ldots, t_N\}$ of $[0, T]$ is arbitrary, then we have
\begin{align*}
\tilde{S}_{0T}(\varphi, \mu) \ge S_{0T}(\varphi, \mu).
\end{align*}
The equality follows from the fact that equation \eqref{Eq15} becomes an equality when $(\varphi, \mu)$ is piecewise linear w.r.t. the  partition $\{t_0, t_1, \ldots, t_N\}$. This completes the proof. \qed

In what follows, we assume that $\epsilon^{-1}\hat{\delta} \rightarrow 0$ as $\hat{\delta} \rightarrow 0$. Then, let us drop $\hat{\delta}$ in the notation and consider a variation $\mathbb{Q}_{z_0, k}^{\epsilon}$ of $\mathbb{P}_{z_0, k}^{\epsilon}$ that is governed by the same initial value and evolution except that its transmutation intensity $c_{km}(t)$ depends on time rather than state position. Then, the corresponding Legendre transform for $\mathbb{Q}_{z_0, k}^{\epsilon}$ in $(p, \alpha)$ is given by
\begin{align*}
\hat{\eta}(t, z, q, \beta) =\sup_{p \in \mathbb{R}^{(n+m)},\, \alpha \in \mathbb{R}^K}  \bigl[ q \cdot p + \beta \cdot \alpha - \hat{\lambda}(t, z, p, \alpha) \bigr], \,\,\, z, q \in \mathbb{R}^{(n+m)}, \,\,\, \beta \in \mathbb{R}^K,
\end{align*}
where the principal eigenvalue $\hat{\lambda}(t, z, p, \alpha)$ is associated with the following matrix $\bigl(\hat{H}_{km}(t, z, p,\alpha)\bigr)$ with
\begin{align*}
\hat{H}_{km}(t, z,p,\alpha) = \left \{ \begin{array}{l}
  c_{km}(t), \hspace{2.52in} \text{if} \quad m \neq k,\\
  \bigl[p \cdot \tilde{a}_k(z) p/2 + p \cdot \tilde{F}_k(z) + \alpha_k \bigr]  - \sum\nolimits_{j: j \neq k} c_{kj}(t), \quad \text{if} \quad m = k.
\end{array} \right.
\end{align*}
where $\tilde{a}_k(z) = \tilde{\sigma}_k(z)\,\tilde{\sigma}_k^T(z)$ for $k =1, 2, \ldots, K$.

Let $T > 0$ be fixed and define
\begin{align*}
\hat{S}_{0T}(\varphi, \mu) = \left \{ \begin{array}{l}
  \int_0^T \hat{\eta}\bigl(s, \varphi(s), \dot{\varphi}(s), \dot{\mu}(s)\bigr) ds, \,\,\, \text{if} \,\, \varphi \in C(\mathbb{R}^{(n+m)}) \,\, \text{and} \,\, \mu \in C_{+}(\mathbb{R}^K)\,\, \text{a.c.}\\
 +\infty \hspace{1.65 in}\text{otherwise}.
\end{array} \right.
\end{align*}

Let $\epsilon^{-1}\hat{\delta} \rightarrow 0$ as $\hat{\delta} \rightarrow 0$, then we have the following result.
\begin{proposition} \label{P4}
The action functional for the family of processes $(z_t^{\epsilon}, r_t^{\epsilon})$ w.r.t. $\mathbb{Q}_{z_0, k}^{\epsilon}$ as $\epsilon \rightarrow 0$ is ${\epsilon}^{-1}\hat{S}_{0T}(\varphi, \mu)$ in the uniform topology.
\end{proposition}

{\em Proof}: 
Note that when constant diffusivity and drift coefficients, this lemma can be checked by following the same argument as that of Proposition~\ref{P3}. Here, the standard argument is to freeze the diffusivity and drift coefficients for smaller and smaller durations and then updating them afterwards and extending the statement to the general case for $\mathbb{Q}_{z_0,k}^{\epsilon}$ (cf. \cite[Section~6]{Var84}). Hence, we omit the details.\qed

{\em Proof of Proposition~\ref{P2}}: 
Let us define
\begin{align*}
 I(\epsilon, \delta) = \mathbb{P}_{z_0,k}^{\epsilon} \Bigl \{\max_{0 \le s \le T} \bigl\vert z_s^{\epsilon} - \varphi(s) \bigr\vert < \delta, \max_{0 \le s \le T} \bigl\vert r_s^{\epsilon} - \mu(s) \bigr\vert < \delta\Bigr\}.
\end{align*}
By a standard argument from the theory of large deviations, then it is easy to see that
\begin{align}
\text{if} \quad S_{0T}(\varphi, \mu) < \infty, \quad \text{then} \quad \lim_{\delta \rightarrow 0} \lim_{\epsilon \rightarrow 0}  I(\epsilon, \delta) = - S_{0T}(\varphi, \mu), \label{Eq16}
\end{align}
that is, the exponential tightness condition (e.g., see \cite[Section~3.3]{FreWe84}); and, hence, there exists a sequence of $\Omega_j$ of compact subsets of the trajectories on $[0, T]$ in the uniform topology such that
\begin{align*}
\lim_{j \rightarrow \infty} \lim_{\epsilon \rightarrow 0}  \epsilon \log \mathbb{P}_{z_0,k}^{\epsilon} \Bigl\{ \bigl( z_{\cdot}^{\epsilon}, r_{\cdot}^{\epsilon} \bigr) \notin \Omega_j \Bigr\} = -\infty.
\end{align*}
On the other hand, fix $T > 0$, let $\mathbb{Q}_{z_0, k}^{\epsilon}$ be the variation of $\mathbb{P}_{z_0, k}^{\epsilon}$ associated with $c_{mk}(\varphi(t))$ (i.e., noting that  $\mathbb{Q}_{z_0, k}^{\epsilon}$ is absolutely continuous w.r.t. $\mathbb{P}_{z_0, k}^{\epsilon}$). Then, define the following change of measure in the space of trajectories on $[0, T]$, i.e., the Radon-Nikodym derivative,
\begin{align*}
\frac{d \mathbb{P}_{z_0, k}^{\epsilon}}{d \mathbb{Q}_{z_0, k}^{\epsilon}} (z,r) &= \exp \frac{1}{\epsilon} \biggl\{\int_0^T \Bigl[c_{\nu_s}(\varphi(s)) - c_{\nu_s}(z_s^{\epsilon})\Bigr] ds \\
 & \quad\quad - \sum\nolimits_{k,\, m,\, k \neq m}  \epsilon\,  \sum\nolimits_{\tau \in \tau_{k,m}} \log \frac{c_{km}(\varphi(\tau))}{c_{lm}(z_{\tau}^{\epsilon})}  \biggr \},
\end{align*}
where $c_i = \sum\nolimits_{j \neq i} c_{ij}$ and  $\tau_{k,m}$ is the set of times in $[0, T]$ when $\nu_s$ changes from type (or mode) $k$ to another mode $m$.\footnote{Note that such a change of measure is also anticipated by the relative density exponential waiting times (e.g. see \cite{GikS72}).} Then, from the continuity of $c_{km}(z)$, it follows that there exists $h(\delta)$ such that $h(\delta) \rightarrow 0$ as $\delta \rightarrow 0$; and, moreover, we have
\begin{align*}
\biggl \vert \int_0^T \Bigl[c_{\nu_s}(\varphi(s)) - c_{\nu_s}(z_s^{\epsilon})\Bigr]  ds \biggr \vert \le h(\delta)
\end{align*}
and 
\begin{align*}
\sup_{0 \le s \le T} \biggl \vert  \log \frac{c_{km}(\varphi(\tau))}{c_{km}(z_{\tau}^{\epsilon})}  \biggr \vert \le h(\delta)
\end{align*}
for $\max_{0 \le s \le T} \bigl \vert z_s^{\epsilon} - \varphi(s) \bigr \vert \le \delta$.

Furthermore, let $\bar{\mathbb{E}}_{z_0,k}^{\epsilon}$ be the associated expectation with $\mathbb{Q}_{z_0,k}^{\epsilon}$ and $\bigl\vert \bigcup_{k,\,m, \,k \neq m} \tau_{k,m} \bigr \vert$ denote the total number of changes of types. Note that $\mathbb{Q}_{z_0,k}^{\epsilon}$ is the absolute continuous w.r.t $\mathbb{P}_{z_0,k}^{\epsilon}$, then we have
\begin{align}
 I(\epsilon, \delta) &= \bar{\mathbb{E}}_{z_0,k}^{\epsilon} \biggl  \{\frac{d \mathbb{P}_{z_0, k}^{\epsilon}}{d \mathbb{Q}_{z_0, k}^{\epsilon}} (z,r) \, \biggl \vert \, \max_{0 \le s \le T} \bigl\vert z_s^{\epsilon} - \varphi(s) \bigr\vert < \delta \,\,\, \text{and} \,\,\, \max_{0 \le s \le T} \bigl\vert r_s^{\epsilon} - \mu(s) \bigr\vert < \delta\biggr\} \notag \\
 &\le \exp \Bigl(\frac{h(\delta)}{\epsilon} \Bigr) \bar{\mathbb{E}}_{z_0,k}^{\epsilon} \biggl \{\exp \biggl( h(\delta) \bigl\vert \bigcup\nolimits_{k,\,m, \,k \neq m} \tau_{k,m} \bigr \vert \biggr) \, \biggl \vert \, \max_{0 \le s \le T} \bigl\vert z_s^{\epsilon} - \varphi(s) \bigr\vert < \delta \notag \\
 & \hspace{2.9 in } \text{and} \,\,\, \max_{0 \le s \le T} \bigl\vert r_s^{\epsilon} - \mu(s) \bigr\vert < \delta\biggr\}. \label{Eq17}
\end{align}
Then, by the H\"{o}lder's inequality, we have
\begin{align*}
 I(\epsilon, \delta) &\le \exp\Bigl(\frac{h(\delta)}{\epsilon} \Bigr) \hat{\mathbb{Q}}_{z_0,k}^{\epsilon} \biggl \{ \max_{0 \le s \le T} \bigl\vert r_s^{\epsilon} - \varphi(s) \bigr\vert < \delta,\, \max_{0 \le s \le T} \bigl\vert r_s^{\epsilon} - \mu(s) \bigr\vert < \delta\biggr\}^{1 - 1/q} \\
 & \hspace{2.5 in} \times \bar{\mathbb{E}}_{x,k}^{\epsilon} \biggl \{\exp \biggl(q h(\delta) \bigl\vert \bigcup\nolimits_{k,\,m, \,k \neq m} \tau_{k,m} \bigr \vert \biggr) \biggr\}^{1/q}
\end{align*}
for all $q > 1$.

Denote $\max_{0 \le s \le T} \Bigl\{c_{km} \, \bigl \vert \, k, m \in\{1, 2, \ldots, K\}, \, k \neq m\Bigr\}$ by $\Xi$. Then, it is easy to check that
\begin{align}
\bar{\mathbb{E}}_{z,k}^{\epsilon} \biggl \{\exp \biggl(q h(\delta) \bigl\vert \bigcup\nolimits_{k,\,m, \,k \neq m} \tau_{k,m} \bigr \vert \biggr) \biggr\} \le \exp \biggl(T \frac{(n-1)\Xi}{\epsilon} \bigl(e^{qh(\delta)} - 1\bigr)\biggr). \label{Eq18}
\end{align}
The above equation together with \eqref{Eq17} imply that
\begin{align*}
\limsup_{\delta \rightarrow 0} \lim_{\epsilon \rightarrow 0}  \epsilon \log  I(\epsilon, \delta) \le - \bigl(1 - 1/q \bigr) \bar{S}_{0T} (\varphi, \mu) =  - \bigl(1 - 1/q \bigr) S_{0T} (\varphi, \mu). 
\end{align*}
Then, letting $q$ to $\infty$ gives the desired upper estimate of \eqref{Eq16}. Note that a lower bound estimate is obtained by a similar argument, when Jensen's inequality is taking the place of the H\"{o}lder's inequality in \eqref{Eq17}. This completes the proof.\qed

{\em Proof of Proposition~\ref{P2}}: 
Let $\Omega_{\delta}$ and $\Omega_{2\delta}$ be $\delta$ and $2\delta$-neighborhoods of the compact set $\Omega \subset G$ with boundaries $\partial \Omega_{\delta}$ and $\partial \Omega_{2\delta}$, respectively. Then, the state-trajectories $z_t^{\epsilon}$, starting from any $z \in G$, $k \in \{1,2, \ldots, K\}$, hit $\partial \Omega_{\delta}$ before $\partial G$ with probability close to one as $\epsilon$ is small enough. This follows from Assumption~\ref{AS1}. Hence, taking into account the strong Markov property of the process $(z_t^{\epsilon}, \nu_{t}^{\epsilon}, \mathbb{P}_{z_0,k}^{\epsilon})$, it is sufficient to prove Proposition~\ref{P2} for $z \in \partial \Omega_{\delta}$, $k \in \{1,2, \ldots, K\}$.

Define the following Markov times $\zeta_0 < \tau_1 < \zeta_1 < \cdots < \tau_n < \zeta_n \cdots$ as follows 
\begin{align*}
\left. \begin{array}{l}
 \zeta_0 = \min \bigl\{t > 0 \,\vert \, z_t^{\epsilon} \in\partial \Omega_{2\delta} \bigr\}\\
 \tau_1 = \min \bigl\{t > \zeta_0 \,\vert \, z_t^{\epsilon} \in\partial \Omega_{\delta} \cup \partial G \bigr\}\\
 \zeta_1 = \min \bigl\{t > \tau_1 \,\vert \, z_t^{\epsilon} \in\partial \Omega_{2\delta} \bigr\}\\
  \hspace{0.75 in} \cdots \\
 \tau_{n+1} = \min \bigl\{t > \zeta_n \,\vert \, z_t^{\epsilon} \in\partial \Omega_{\delta} \cup \partial G \bigr\}\\
 \zeta_{n+1} = \min \bigl\{t > \tau_n \,\vert \, z_t^{\epsilon} \in\partial \Omega_{2\delta} \bigr\}\\
  \hspace{0.75 in } \cdots 
\end{array} \right.
\end{align*}
Next, let us define a Markov chain $(z_t^{\epsilon}, \hat{\nu}_{t}^{\epsilon})$ in the phase space $\bigl\{ \Omega_{\delta} \cup \partial G\bigr\} \times \bigl\{1,2, \ldots, K\bigr\}$. Note that the first exit of $z_t^{\epsilon}$ from $G$ occurs, when the component $z_t^{\epsilon}$ of the chain first time belongs to $\partial G$. Then, using the large deviation estimates for the family of processes $(z_t^{\epsilon}, \nu_{t}^{\epsilon}, \mathbb{P}_{z_0,k}^{\epsilon})$ as $\epsilon \rightarrow 0$, we can show, in the standard way (e.g., see \cite[Chapter~4]{FreWe84}), that $z_t^{\epsilon}$ starting from any $z_0 \in \Omega_{\delta}$ and $k \in \{1,2, \ldots, K\}$ reaches $\partial G$ for the first time to a small neighborhood of the point $z_0 \in \partial G$, introduced in Assumption~\ref{AS2}, with probability close to one as the parameters $\epsilon$ and $\delta$ are small enough, which implies the first statement of Proposition~\ref{P2}.

In order to prove the second statement, we use the fact that the extremal of the variational problem
\begin{align*}
\inf \Bigl\{ I_{0T}(\varphi) \, \bigl\vert \, \varphi(0) \in \Omega, \,\, \varphi(T) = \in \partial G, \,\, T > 0 \Bigr\}
\end{align*}
spends in $\delta$-neighborhood $\partial G_{\delta} = \bigl\{z \in G \, \vert \, \rho(z, \partial G) < \delta \bigr\}$ of $\partial G$ a time of order $\delta$ as $\delta \rightarrow 0$. Note that, with probability close to one as $\delta$ is small, the second component $\nu_{t}^{\epsilon}$ has no jumps during this time; and, hence, $z_t^{\epsilon}$ hits the boundary for the value of the second coordinate $\nu_{t}^{\epsilon}$ such that the transition of  $z_t^{\epsilon}$ from $\partial G_{\delta} \setminus \partial G$ to $\partial G$ is easiest transition, when the second component is equal to $k_0$ defined in Assumption~\ref{AS2}. This completes the proof.\qed

Note that the limiting distributions of $(z_{\tau_{G}^{\epsilon}}^{\epsilon}, \nu_{\tau_{G}^{\epsilon}}^{\epsilon})$ as $\epsilon \rightarrow 0$ also determine the limiting behavior for the solutions of the Dirichlet problem in \eqref{Eq5}, where such a connection is established using the following fact.
\begin{proposition} \label{P5}
Suppose Assumptions~\ref{AS1} and \ref{AS2} hold true. Let $\epsilon^{-1}\hat{\delta} \rightarrow 0$ as $\hat{\delta} \rightarrow 0$, then we have
\begin{align}
\lim_{\epsilon \rightarrow 0} \upsilon_k^{\epsilon}(z) = g_{k_0} (\bar{z}_0), \,\, 1 \le k \le K, \label{Eq13}
\end{align}
uniformly in $z \in \Omega \subset G$,  where $\upsilon_k^{\epsilon}(z)$ is the solution for the Dirichlet problem in \eqref{Eq5}.
\end{proposition}

{\em Proof}: 
The proof easily follows from Proposition~\ref{P2} and the representation
\begin{align*}
\upsilon_k^{\epsilon}(z) = \mathbb{E}_{z_0,k}^{\epsilon} \Bigl\{ \upsilon_{\nu_{\tau_G^{\epsilon}}^{\epsilon}}^{\epsilon}(x_{\tau_G^{\epsilon}}^{\epsilon}) \Bigr\}
\end{align*}
uniformly in $z \in \Omega \subset G$ (cf. \cite[Theorem~3]{EizF90}). Furthermore, note that
\begin{align*}
\mathbb{P}_{z_0,k}^{\epsilon} \Bigl\{ \tau_G^{\epsilon} < \infty \Bigr\} = 1, \,\,\, \text{for any} \,\, z_0 \in \Omega, \,\ k \in \{1,2, \ldots, K\}.
\end{align*}
Thus, taking into account the boundedness (cf. Lemma~\ref{L1}), we have
\begin{align*}
\lim_{\epsilon \rightarrow 0} \upsilon_k^{\epsilon}(z) &= \lim_{\epsilon \rightarrow 0}\mathbb{E}_{z_0,k}^{\epsilon} \Bigl\{ \upsilon_{\nu_{\tau_G^{\epsilon}}^{\epsilon}}^{\epsilon}(z_{\tau_G^{\epsilon}}^{\epsilon}) \Bigr\}\\
&= g_{k_0} (\bar{z}_0) 
\end{align*}
for $k \in \{1, 2, \ldots, K\}$. This completes the proof. 
\qed



\end{document}